\documentclass[11pt]{article}
\newcommand{\paper}[1]{#1}\newcommand{\conference}[1]{}

\usepackage{amssymb}
\usepackage{graphicx}

\newcommand{\R}{{\mathbb R}}
\newcommand{\ip}[2]{\langle #1, #2 \rangle}
\newcommand{\ipT}[2]{\langle #1, #2 \rangle_T}
\newcommand{\abs}[1]{\left\vert #1 \right\vert}
\newcommand{\norm}[1]{\left\|#1\right\|}
\newcommand{\normT}[1]{\left\|#1\right\|_T}
\newcommand{\angT}{\theta _T}
\newcommand{\beqn}{\begin{eqnarray*}}
\newcommand{\eeqn}{\end{eqnarray*}}
\newcommand{\text}[1]{\hbox{\rm \ #1\ \/}}
\newcommand{\be}[1]{\begin{equation}\label{#1}}
\newcommand{\ee}{\end{equation}}
\newcommand{\RE}{\mbox{Re\,}}
\newcommand{\gi}{\gamma _\infty }
\newcommand{\xstar}{x^*}

\title{A notion of passivity gain and a generalization of the ``secant condition'' for stability} 
\author{Eduardo D. Sontag\footnote{Email: sontag@math.rutgers.edu}\\
Dept.\  of Mathematics, Rutgers University, New Brunswick, NJ}

\begin{document}

\maketitle

\begin{abstract}
A generalization of the classical secant condition for the stability of
cascades of scalar linear systems is provided for passive systems.  The key is
the introduction of a quantity that combines gain and phase information for
each system in the cascade.  For linear one-dimensional systems, the known
result is recovered exactly.
\end{abstract}

\section{Introduction}

An often-used tool in the analysis of biological feedback loops is the 
{\em secant condition\/} for linear stability; see the classical
papers by Tyson and Othmer \cite{Tyson-Othmer} and Thron \cite{Thron}, as well
as the recent paper~\cite{cauchy}.
Consider a matrix of the following form:
\[
\pmatrix{-\alpha _1&0&\ldots &0&-\beta _1\cr
      \beta _2&-\alpha _2&\ldots &0&0\cr
      \vdots&\vdots&  & &\vdots\cr
      0&0&\ldots &\beta _n&-\alpha _n}
\]
where all $\alpha _i>0$ and all $\beta _i>0$.
Then, the secant condition states that the matrix is Hurwitz provided that:
\[
{\beta _1\ldots \beta _n\over\alpha _1\ldots \alpha _n}
   < \left(\sec {\pi \over n}\right)^n \,.
\]
In essence, this says that a stable system with distinct real eigenvalues and
no zeros tolerates negative feedback with a gain much larger than that
provided by the small-gain theorem: the corresponding small-gain estimate
would just have a ``1'' in the right-hand side.  (The secant expression is
always bigger than one.  It is singular at $n=2$ --which it should be, since
then the matrix is always Hurwitz-- and it equals 8 for $n=3$, 4 for $n=4$, and
$\approx2.88$ for $n=5$, and tends monotonically to 1 as $n\rightarrow \infty $.  The bound
is achieved exactly when all the $\alpha _i$'s are the same.)  The condition is
useful because certain models of biological systems are stable for gains
larger than those predicted by a simple application of the small-gain theorem.
The secant takes advantage simultaneously of phase and gain information on the
open-loop system.

We provide here a generalization of the secant condition to cascades of output
strictly passive (OSP) systems.  We do so in such a manner that, in the
special case in which each system is linear and one-dimensional, the classical
result is recovered.  (For simplicity, we restrict ourselves to systems with
scalar inputs and outputs, but it is obvious how to generalize to more
arbitrary systems, as long as inputs and outputs have the same dimension.)

The generalization is based on systematic use of a ``gain'' associated to OSP
systems.  It would seem that the use of this quantity might be useful for many
other problems as well.

This note is organized as follows.  Section~\ref{intro-sec} introduces the
basic concepts and states the main result, which is proved in
Section~\ref{proof-sec} (the proof is actually very easy, given the
definitions).
Section~\ref{extensions-sec} briefly mentions some extensions of the basic
formalism, and Section~\ref{remarks-sec} collects several facts concerning
secant gains for the special case of linear systems.

\section{Notations, Definitions, and Statement of Main Result}
\label{intro-sec}

As usual, the extended space $L^2_e(0,\infty )$ denotes the set of signals (thought
of as time functions) $w:[0,\infty )\rightarrow \R$ which have the property that
each restriction  $w_T = w|_{[0,T]}$ is in $L^2(0,T)$, for every $T>0$.
Given an element $w\in L^2_e(0,\infty )$ and any fixed $T>0$, one writes
$\normT{w}$ for the $L^2$ the norm of this restriction $w_T$, and given
two functions $v,w\in L^2_e(0,\infty )$ and any fixed $T>0$,
the inner product of $v_T$ and $w_T$ is denoted by $\ipT{v}{w}$ .
In any Hilbert space, one defines the angle $\theta (v,w)\in [0,\pi ]$
between two elements $v,w$ by the formula
\[
\cos \theta (v,w) = \frac{\ip{v}{w}}{\norm{v}\norm{w}}
\]
if $v$ and $w$ are nonzero, and zero otherwise.
Given $v,w\in L^2_e(0,\infty )$ and any fixed $T>0$, we will write
$\angT(v,w)$ instead of $\theta (v_T,w_T)$, to denote the angle between the
restrictions of the signals to $[0,T]$.

We consider continuous-time finite-dimensional systems $\dot x=f(x,u)$, $y=h(x)$
in the usual sense of control theory (e.g.~\cite{mct}), with scalar valued
inputs and outputs, and state space $\R^n$, and assume always that the system
is $L^2$-well-posed, in the sense that for each $u\in L^2_e(0,\infty )$ and initial
state $x(0)=0$ there is a unique solution $x(\cdot )$ defined for all $t\geq 0$
and the corresponding output $y(t)=h(x(t))$ is also in $L^2_e(0,\infty )$.
We call $(u,y)$ is an {\em input/output (i/o) pair\/} of the system.

We recall the standard notion of an {\em output strictly passive\/} (``OSP''
for short) system, as given in textbooks such as \cite{khalil,vdS,Vidyasagar}.
A system is OSP if there is some $\gamma >0$ such that,
for every i/o pair $(u,y)$,
\be{OSP}
\normT{y}^2 \;\leq \; \gamma  \, \ipT{u}{y}
\ee
for all $T>0$.  (Allowing an additive constant in the inequality is useful
when dealing with arbitrary initial states.  As we will study zero-state
responses, we do not include a constant.)

If a system is OSP, we call the smallest $\gamma $ as in~(\ref{OSP})
the {\em secant gain\/} of the system, and denote it as $\gamma _s$.
(There is a smallest such $\gamma $, since the set of $\gamma $'s that
satisfy~(\ref{OSP}) is a closed set.)

An equivalent definition of $\gamma _s$ is as the smallest $\gamma $ with the property
that 
\[
\normT{y}^2\leq \gamma \normT{u}\normT{y} \angT(u,y),
\]
 or equivalently:
\be{OSP-angle}
\normT{y} \;\leq \; \gamma  \,\normT{u} \,\cos\angT(u,y)
\ee
for all $T>0$ and all i/o pairs.
Since~(\ref{OSP}) implies that $\ipT{u}{y}\geq 0$ for all i/o pairs
and all $T$, for OSP systems we always think of the angle as lying
in the interval $[0,\pi /2]$, and the cosine is nonnegative.

The Cauchy-Schwartz inequality applied to~(\ref{OSP}) gives
$\normT{y}\leq \gamma \normT{u}\leq \gamma \norm{u}$ for all $T>0$,
so in particular $y\in L^2$ if $u\in L^2$,
and an OSP system necessarily has finite $L^2$-induced (or ``$H_\infty $'') gain
$\gi\leq \gamma _s$ (we remark later that this inequality is in general a strict one).
Just as the $L^2$ gain is the supremum of the expressions
$\normT{y}/\normT{u}$ over all $T$ and all i/o pairs with nonzero $u$, the
secant gain is obtained by maximizing $\sec\angT(u,y)\normT{y}/\normT{u}$,
hence our terminology.

If $u\in L^2$, so that also $y\in L^2$, taking limits in~(\ref{OSP}) gives
\be{OSP-L2}
\norm{y}^2\;\leq \;\gamma \ip{u}{y}\,.
\ee
Conversely, if $u\in L^2\Rightarrow y\in L^2$ and~(\ref{OSP-L2}) is true for all $u\in L^2$,
then~(\ref{OSP}) holds.
This is a routine exercise in causality, as follows.
Pick any i/o pair $(u,y)$ and any $T>0$.  Let $v\in L^2$ be input which equals
$u$ on $[0,T]$ and is zero for $t>T$, and $z$ the output corresponding to $v$.
Since $v\in L^2$, also $z\in L^2$.  By causality, $z$ restricted to $[0,T]$ is the
same as $y$ restricted to $[0,T]$, so $\ip{v}{z}=\ipT{u}{y}$, and
$\normT{y}=\normT{z}$. 
Therefore $\normT{y}^2=\normT{z}^2\leq \norm{z}^2\leq \gamma \ip{v}{z}=\gamma \ipT{u}{y}$,
and indeed~(\ref{OSP}) is verified.

\begin{figure}[h,t]
\begin{center}
\paper{\setlength{\unitlength}{2000sp}}
\conference{\setlength{\unitlength}{1500sp}}
\begin{picture}(9381,1011)(139,-3160)
\put(901,-2761){\circle{300}}
\put(1801,-3061){\framebox(900,600){}}
\put(3301,-3061){\framebox(900,600){}}
\put(4801,-3061){\framebox(900,600){}}
\put(8101,-3061){\framebox(900,600){}}
\put(2701,-2761){\vector( 1, 0){600}}
\put(4201,-2761){\vector( 1, 0){600}}
\multiput(5701,-2761)(117.07317,0.00000){21}{\line( 1, 0){ 58.537}}
\put(8101,-2761){\vector( 1, 0){0}}
\put(9001,-2761){\line( 1, 0){300}}
\put(9301,-2761){\line( 0, 1){600}}
\put(9301,-2161){\line(-1, 0){8400}}
\put(901,-2161){\vector( 0,-1){450}}
\put(1051,-2761){\vector( 1, 0){750}}
\put(151,-2761){\vector( 1, 0){600}}
\put(1051,-2461){$-$}
\put(301,-3061){$u$}
\put(2851,-3061){$y_1$}
\put(4351,-3061){$y_2$}
\put(9151,-3061){$y_n$}
\put(451,-2611){$+$}
\end{picture}
\caption{Closed-loop system}
\label{fig-system}
\end{center}
\end{figure}
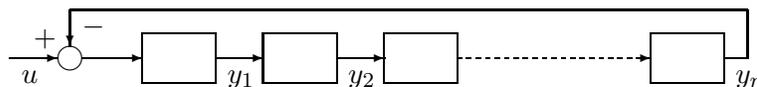

We wish to analyze the stability of the closed-loop system $\dot x=f(x,u-h(x))$
obtained under negative unity feedback.
Specifically, we study a cascade of $n$ subsystems, as shown in the diagram in
Figure~\ref{fig-system}
and subject to unity negative feedback.
Such cascades appear frequently in control theory as well as in biological
applications, and, when components are one-dimensional, tend to have
especially good dynamical properties such as the validity of the
Poincar\'e-Bendixson Theorem (\cite{cyclic}).
We will assume that the $i$-th system has a secant gain $\gamma _i$, and we write
$y_i$ for the output of the $i$th subsystem.
We also assume well-posedness of the closed-loop.


The main result is as follows:

\noindent
{\bf Theorem.}
Suppose that 
\[
\gamma _1\gamma _2\ldots \gamma _n \;<\; \left(\sec {\pi \over n}\right)^n \,.
\]
Then the cascade is $L^2$-stable: there is a number $c$ so that
\[
\normT{y_n} \;\leq \; c \normT{u}
\]
for all input/output pairs in the cascade and all $T>0$.

\medskip

Of course, this property implies as well that every $\normT{y_i}$ is bounded
by some linear function of $\normT{u}$, and that the signals $y_i$ belong to
$L^2$ if $u\in L^2$.  

For the special cases $n=1$ and $n=2$ (secant is infinite), we interpret the
inequality in the theorem as saying that the condition holds for any possible
values of the $\gamma _i$'s.  For $n=2$, therefore, the theorem is simply a
restatement of the Passivity Theorem as given e.g.\ 
in~\cite{vdS}, Theorem~2.2.15, Part a (using only the input $u$).
The Passivity Theorem also includes usually a statement (``Part b'' in the
citation) regarding the case in which the first system is OSP and the second
one is only passive, meaning that only $\ipT{u}{y}\geq 0$ is known for all i/o
pairs.  We comment later on this fact.

The assumption that the initial state of the cascade is $x(0)=0$ is easy to
dispose of, assuming appropriate reachability of the cascade, as routinely
done in going from input/output stability to state space stability, and
Barb\u{a}lat's Lemma combined with either reachability 
or detectability arguments can be used to show convergence of internal states
to zero.  As an illustration, we state just one such corollary:

\noindent
{\bf Corollary.}
Suppose that the condition in the Theorem is verified, that the composite
system shown in Figure~\ref{fig-system} is zero-reachable and that each
subsystem is input to state $L^2$-stable.  Then the system with no inputs
($u=0$) has the property that all solutions converge to $x=0$. 

\section{Extensions}
\label{extensions-sec}

We have formulated the results in terms of state-space systems only in order
to be concrete.  One could equally well consider arbitrary operators
$L^2\rightarrow L^2$, or even just relations $R$ on $L^2\times L^2$, where an ``i/o pair'' is
by definition any element of $R$, and define secant gain $\gamma _s$ as the
smallest number so that~(\ref{OSP}) holds for all $T$ and all i/o pairs.  Nor
is it needed for the inputs and outputs to be scalar-valued; one may consider
values on arbitrary Hilbert spaces, with inner product and norms taken
pointwise in that space.  
More generally, functions of time are not required: one could consider an
arbitrary Hilbert space $H$ and simply ask that $u$ and $y$ belong to $H$.
(To be precise, one needs a Hilbert space together with a resolution of the
identity, in order to be able to be able to have a concept of ``restriction''
of $u$ and $y$ to subintervals; this is the formalism of 
{\em resolution spaces\/} developed in~\cite{saeks}.) 
Even more generally, if one has a system in which inputs $u$ and outputs $y$
are known to lie in a specific subset $S\subseteq H$, then $\gamma _s$ can be defined in
terms only of i/o pairs that lie in $S$; the validity of the main theorem is
not affected, since it is just an algebraic statement about norms and inner
products. 

Let us discuss a simple example of an operator defined only on subsets, 
which is of interest in biomolecular applications
(``Michaelis-Menten kinetics'').
Suppose that $S$ is the set of all $L^2$ maps
$w:[0,\infty )\rightarrow [-a,\infty )$ with any fixed $a>0$, and that we consider the
function $\ell:[-a,\infty )\rightarrow \R$ given by
\be{mapping-eg}
\ell(r) = \frac{Vr}{K+a+r}
\ee
(with $K,V>0$ some constants)
and the operator $u\mapsto y$ defined on $S$, where $y(t)=F(u)(t) = \ell(u(t))$.
This is an example of a ``sector'' nonlinearity.  The analysis of sector
nonlinearities is routine in passivity theory.
The operator $F$ is OSP and has $\gamma _s=V/K$, because we have, for all
$r\in [-a,\infty )$: 
\beqn
[\ell(r)]^2 
  &=& \frac{V}{K+a+r} \frac{Vr^2}{K+a+r}\\
  &\leq & \frac{V}{K} \frac{Vr^2}{K+a+r}
  \;=\; \frac{V}{K} r \ell(r)
\eeqn
(since $K+a+r\geq K$), and thus
\beqn
\normT{y}^2 &=&
\int_0^T\ell(u(t))^2 \,dt\\
&\leq &\frac{V}{K}\int_0^T u(t)\ell(u(t))\,dt
\;=\;\frac{V}{K}\ipT{u}{y}
\eeqn
so $\gamma _s\leq V/K$, and the equality is verified when $u(t)\equiv -a$.

Stability in the $L^2$ sense is only appropriate when dealing with equilibria
associated to zero signals.  However, the framework described here can be
easily extended to more general situations.  These extensions are of interest,
particularly, when dealing with problems in biology and chemistry, where
quantities represent concentrations of substances, and hence are always
nonnegative.  We now describe briefly how this extension can be accomplished.

Suppose that one wishes to study a system
\beqn
\dot x &=& f(x,u)\\
y&=&h(x)
\eeqn
under the feedback law $u=-y$, and that there is a steady state $\xstar$
for this closed-loop system:
\[
f(\xstar,-h(\xstar))=0
\]
whose stability is of interest to analyze.
We assume that the states $x(t)$ evolve in some subset $S$ of $\R^n$, for
example the positive orthant 
$\R_+^n = \{(x_1,\ldots ,x_n),\, x_i\geq 0\,\forall i\}$,
and inputs $u$ of the open-loop system take values on some set $U$.
(In order  for the closed-loop system to make sense, one should then have that
$-h(S)\subseteq U$, of course.)
We perform a change of variables
$
z = x-\xstar
$
and define the new system
\beqn
\dot z &=& g(z,v) \;=\; f(z+\xstar,v-h(\xstar))\\
w&=&\ell(z) \;=\; h(z+\xstar)-h(\xstar),
\eeqn
with states $z(t)$ in the state-space $\{x-\xstar,\, x\in S\}$,
inputs $v(t)$ in the input-value space $\{u+h(\xstar),\,u\in U\}$,
and outputs $w(t)$.
Note that $g(0,0)=0$.  
Applying the feedback $v=-\ell(z)$ results in
\[
\dot z = g(z,-\ell(z)) = f(z+\xstar,-h(z+\xstar)).
\]
Therefore, for each solution $x(t)$ of $\dot x = f(x,-h(x))$, the vector function
$z(t)=x(t)-\xstar$ satisfies $\dot z = g(z,-\ell(z))$, and conversely, each
solution of the latter system arises from the former.
Proving that solutions of $\dot x = f(x,-h(x))$ converge to $\xstar$ is then
equivalent to proving that the solutions of the new system converge to the
equilibrium $z=0$.  Thus we have reduced the analysis to the case treated in
this paper.

For example, suppose that we wish to study a positive system, that is,
a system whose state state space is $\R^n_+$ and inputs are
also nonnegative.
Furthermore, suppose that, as is often the case in biological feedback loops,
one wishes to study an inhibitory feedback of the form
\[
u = \frac{M}{K+x_n}
\]
where $M$ and $K$ are some positive constants and $x_n$ is the $n$th
coordinate of the state, that is to say, we have $h(x)=-M/(K+x_n)$.
In terms of the variables $z$, we have the output
\beqn
w &=& \ell(z) = h(z+\xstar)-h(\xstar)\\
  &=&\frac{M}{K+\xstar_n}-\frac{M}{K+(z_n+\xstar_n)}
  = \frac{Vz_n}{K+\xstar_n+z_n}
\eeqn
which is the function in~(\ref{mapping-eg}) with $a=\xstar$
and $V = M/(K+\xstar_n)$.
Since $x_n(t)$ is nonnegative, the state variable $z_n(t)$ takes values
in $[-\xstar,\infty )$.
Thus, we may view the closed-loop system as built from cascading the original
system (which may itself be a cascade of several subsystems) with the static
system ``$y=\ell(u)$'', which has $\gamma _s=V/K$, and the previous analysis
applies. 

This is all particularly simple for a linear system $\dot x=f(x,u)=Ax+Bu$.
Positivity amounts to asking that all the off-diagonal entries of $A$ as well
as all entries of $B$ are nonnegative 
(see e.g.~\cite{patrick,farina-rinaldi-book}).
Since the system is linear and $A\xstar-Bh(\xstar)=0$, we have that
$g(z,v) = A(z+\xstar)+B(v-h(\xstar)) = Az+Bv$, so the same open loop system
results, except that now we are interested in the stability of $z=0$.

\section{Proof of Main Result}
\label{proof-sec}

Given an external input $u$, the solutions of the closed-loop system with
initial state zero are so that the signals $y_i$ have the following properties:
\beqn
\normT{y_1}^2 &\leq & \gamma _1 \, \ipT{u+y_0}{y_1}\\
\normT{y_2}^2 &\leq & \gamma _2 \,\ipT{y_1}{y_2}\\
\vdots\\
\normT{y_n}^2 &\leq & \gamma _n \,\ipT{y_{n-1}}{y_n}
\eeqn
for every $T>0$, where we are writing $y_0=-y_n$.
We expand $\ipT{u+y_0}{y_1}=\ipT{u}{y_1}+\ipT{y_0}{y_1}$, and use the
Cauchy-Schwartz inequality for the first term, upper-bounding it by
$\normT{u}\normT{y_1}$.  
Replacing now each $\ipT{y_{i-1}}{y_i}$ by 
$\normT{y_{i-1}}\normT{y_i}\cos\angT(y_{i-1},y_i)$ and dividing by
$\normT{y_i}$ (assumed nonzero; otherwise, there will be nothing to prove),
we have these estimates:
\beqn
\normT{y_1} &\leq & \gamma _1 \,\normT{y_0}\cos\angT(y_0,y_1) 
                   \,+\, \gamma _1\,\normT{u}\\
\normT{y_2} &\leq & \gamma _2\, \normT{y_1}\cos\angT(y_1,y_2) \\
\vdots\\
\normT{y_n} &\leq & \gamma _n\, \normT{y_{n-1}}\cos\angT(y_{n-1},y_n) 
\eeqn
from which we conclude, by recursively substituting the estimates starting
from the last one backward towards the first, that:
\[
\normT{y_n} \;\leq \; \kappa  \normT{y_n} \, +\, \alpha  \normT{u}
\]
where
\[
\alpha \, =\, \gamma _1\gamma _2 \ldots \gamma _n\cos\angT(y_1,y_2) \ldots .\cos\angT(y_{n-1},y_n)
\]
and
\[
\kappa  \;=\; \alpha  \, \cos\angT(y_0,y_1)\,.
\]
It is enough to show that $\kappa <1$, since then we can write 
$(1-\kappa )\normT{y_n}\;\leq \;\alpha \normT{u}$, and therefore the result holds with
$c=\alpha /(1-\kappa )$.  
Let us fix $T$ and write $\theta _i:=\angT(y_{i-1},y_i)$ for $i=1,\ldots ,n$.
We must show, then, that
\be{main-ineq}
\cos \theta _1\ldots \cos\theta _n \leq  \left(\cos {\pi \over n}\right)^n \,.
\ee
The angles $\theta _i$ all lie in $[0,\pi /2]$, for each $i=2,\ldots ,n$, since each
system is OSP; thus $\cos\theta _i\geq 0$ for all such $i$.  However, it is possible
that $\cos\theta _1<0$, since all that is known is that $\ipT{u+y_0}{y_1}\geq 0$,
not that $\ipT{y_0}{y_1}\geq 0$.
But if $\cos\theta _1<0$, then~(\ref{main-ineq}) is true because the left-hand side
is $\leq 0$ and the right-hand side is positive.  So, in order to
prove~(\ref{main-ineq}), we may assume from now on that all $\theta _i\in [0,\pi /2]$.

We prove, more generally, this fact about Hilbert spaces: suppose given vectors
$v_0,v_1,\ldots ,v_n$ such that $\ip{v_i}{v_{i+1}}\geq 0$, and $v_0=-v_n$.
Let $\theta _i\in [0,\pi /2]$ be the angle between $v_{i-1}$ and $v_i$.
Then~(\ref{main-ineq}) holds.  
Intuitively, the property that the start and end vector are
at angle $\pi $ means that the consecutive vectors cannot be too close in
angle, and therefore at least some of the angles must be large, and hence have
small cosine, and the largest possible value is achieved when all angles are
the same.

To prove this general fact, without loss of generality, we may assume
that all the $v_i$ are unit vectors (since only angles matter).
Notice that $\sum_i\theta _i\geq \pi $.
This is because, for any three unit vectors, $\theta (u,v)+\theta (v,w)\geq \theta (u,w)$,
since we can view the angle as the geodesic distance in a sphere, and
apply the triangle inequality; inductively applied starting from $v_0$, we get
that $\sum_i\theta _i\geq \theta (v_0,v_n)=\pi $.
Now, we have also this algebraic fact:
\[
\cos\theta _1\ldots \cos\theta _n \leq  \left(\cos\frac{\theta _1+\ldots +\theta _n}{n}\right)^n
\]
which follows by noticing that the function $f(x)=-\ln\cos x$ is convex 
for $x\in [0,\pi /2)$, applying Jensen's inequality to obtain
$f(\sum_i \theta _i/n) \leq  (1/n) \sum_i f(\theta _i)$, and taking exponentials.
Together with
$\sum_i\theta _i\geq \pi $, using that $\pi /n\leq (\theta _1+\ldots +\theta _n)/n\leq \pi /2<\pi $
(recall that each $\theta _i\in [0,\pi /2]$), and using that $\cos$ decreases on
$[0,\pi ]$, we conclude:
\[
\left(\cos\frac{\theta _1+\ldots +\theta _n}{n}\right)^n \leq  \left(\cos {\pi \over n}\right)^n
\,.
\]
This completes the proof of the Theorem.

To prove the Corollary, we provide a standard argument, as done e.g.\ 
in~\cite{mct}, Theorem~33.  Pick any initial state $x_0$ and consider the
solution $x(\cdot )$ of the closed-loop system $\dot x=f(x,u-h(x))$ with input $0$
and $x(0)=x_0$. 
Zero-reachability means that there is some finite-time input $u_0:[0,T]\rightarrow \R$
such that, if $z_0(\cdot )$ solves the closed-loop equations $\dot z=f(z,u-h(z))$ with
initial state $z_0(0)=0$ and this input $u_0$ on the interval $[0,T]$, then
$z_0(T)=x_0$. 
Consider now the input $u$ obtained by the formula $u(t)=u_0(t)$ for $t\leq T$
and $u(t)\equiv 0$ for $t>T$, and let $z(\cdot )$ be the solution with initial state
$z(0)=0$ and this input $u$; by causality, $z(t)=z_0(t)$ for $t\leq T$, and hence
$z(T)=x_0=x(0)$, from which it follows that $z(t+T)=x(t)$ for all $t\geq 0$.
Showing $x(t)\rightarrow 0$ as $t\rightarrow \infty $ is the same as showing $z(t)\rightarrow 0$ as $t\rightarrow \infty $.
Let $y_i$ be the outputs of the subsystems when using input $u$ (and zero
initial state).  Since $u\in L^2$ and $\norm{y}\leq c\norm{u}<\infty $,
we have that $y_i\in L^2$ for each of the intermediate outputs.
Since each subsystem is input to state $L^2$-stable, meaning that $L^2$ inputs
(and zero initial state) produces $L^2$ state trajectories, we have that the
complete state $z$ is in $L^2$.
Finally, as $z$ is a trajectory of a semiflow in finite dimensions, we must
have that $z(t)\rightarrow 0$, by a Barb\u{a}lat's Lemma type of argument (see
e.g.~\cite{meagre}). 

Finally, we review in the present context a weaker version that applies when
$n=2$, basically part of the statement of the classical Passivity Theorem.  
Suppose that the first system is OSP but the second system is only known to be
passive, in the sense that no estimate $\normT{y_2}^2\leq  \gamma _2 \ipT{y_1}{y_2}$
may hold, but we do know that $\ipT{y_1}{y_2}\geq 0$ for all $T>0$.
Then, $y_0=-y_2$ implies that:
\beqn
\normT{y_1}^2 &\leq & \gamma _1 \ipT{u+y_0}{y_1}\\
&=& \gamma _1 \ipT{u}{y_1} - \gamma _1\ipT{y_2}{y_1} \;\leq \; \gamma _1 \ipT{u}{y_1}
\eeqn
and so the system {\em with output\/} $y_1$ is OSP, and in particular,
$L^2$ stable.  If, in addition, the second system is also $L^2$ stable, then
stability to $y_2$ holds as well.

\section{Linear Systems}
\label{remarks-sec}

The condition that a system be OSP is of course a restrictive one, but the
concept of OSP system is thoroughly well-studied, and examples of passive
systems abound, especially, but not only, for linear systems.  
We collect here some facts, mostly well-known, regarding the linear case.

For a stable linear system with transfer function $G(s)$, the secant gain can
be characterized as the smallest $\gamma $ such that
\be{OSP-linear}
\abs{G(i\omega )}^2\;\leq \;\gamma \,\RE G(i\omega ) \quad\forall\, \omega \in \R \,.
\ee
A proof is as follows.
First of all, squaring the expression below and expanding 
$\ipT{y-(\gamma /2)}{y-(\gamma /2)}$, one easily sees that the definition of OSP system
is equivalent to the requirement that
\be{OSP-2}
\normT{y-(\gamma /2)u} \;\leq \; (\gamma /2)\normT{u}
\ee
for all i/o pairs and all $T$, which means $\gamma _s$ is the smallest number such
that the $L^2$-induced norm of $u\mapsto y-(\gamma /2)u$ is $\leq \gamma /2$.
For linear systems, induced $L^2$-induced norm corresponds to $H_\infty $ gain,
that is to say, $\gamma _s$ is the smallest number so that 
$\sup_{\omega \in \R}\abs{G(i\omega )-(\gamma /2)}\leq \gamma /2$.
Writing $\abs{G(i\omega )-(\gamma /2)}^2=(G(i\omega )-(\gamma /2))(\overline{G(i\omega )}-(\gamma /2))$
and expanding, one has~(\ref{OSP-linear}).

An equivalent formulation of~(\ref{OSP-linear}) is via the following
analog of the estimate~(\ref{OSP-angle}):
\be{OSP-linear-angle}
\abs{G(i\omega )}\;\leq \;\gamma \,\cos \theta (G(i\omega )) \quad\forall\, \omega \in \R \,.
\ee
where we are denoting now by $\theta (\mu )$ the argument of a complex number $\mu $.
Since $G$ is analytic on $\RE \lambda \geq 0$ (stability), the maximum modulus
principle for analytic functions implies that same estimate is obtained when
maximizing not merely over $\lambda =i\omega $ purely imaginary, but also over all
complex numbers with nonnegative real part.  

If we write $G(s)=p(s)/q(s)$ as a quotient of two polynomials,
condition~(\ref{OSP-linear}) can be also written as 
\[
\abs{p(i\omega )}^2\;\leq \;\gamma \,\RE [p(i\omega )\overline{q(i\omega )}]\,.
\]
For example, for a one-dimensional system $\dot x=-\alpha x+\beta u$ with output $y=x$,
the transfer function is $\beta /(s+\alpha )$, so that $p(i\omega )=\beta $ and
$\RE[p(i\omega )\overline{q(i\omega )}]=\alpha \beta $ for any $\omega $, from which it follows that
$\gamma _s=\beta /\alpha $, and the classical result is obtained.
On the other hand, as is well-known for OSP systems, $G(s)$ must have relative
degree at most one (the condition $\RE G(i\omega )\geq 0$ is otherwise violated).
Therefore, cascades, as studied here, of two or more such one-dimensional
systems are not OSP themselves.

For linear systems, a sufficient condition for a system to be OSP is that its
transfer function $G(s)$ be strictly positive real (SPR), meaning that
$G(s-\varepsilon )$ is positive real for some $\varepsilon >0$, or equivalently
(see e.g.~\cite{khalil}, Lemma 10.1) that it be stable (all poles have
negative real part) and satisfy $\RE G(i\omega ) >0$ for all $\omega \in \R$
and $\lim_{\omega \rightarrow \infty }\omega ^2\RE G(i\omega ) > 0$.
(Note that our transfer functions are strictly proper, by definition, since we
are considering state-space systems with no direct i/o term; for non-strictly
proper transfer functions, the condition is slightly different.)
This provides a large class of examples; for instance, any transfer function
of the form $(s+\alpha )/(s^2+as+b)$ with $b>0$ and $0<a<2\sqrt{b}$ is SPR if and
only if $0<\alpha <a$ (\cite{khalil}, Exercise 10.1).
That SPR implies OSP can be proved using the
Kalman-Yakubovich-Popov (KYP) Lemma.  The converse implication does not hold:
$s/(s^2+s+1)$ is not SPR, since it fails the test just quoted with
($a=b=1$, $\alpha =1/2$) or just by noting that there is an imaginary axis zero,
since $\RE G(0)=0$,
but it is OSP, since $\abs{p(i\omega )}/\RE [p(i\omega )\overline{q(i\omega )}]\equiv 1<\infty $.

More generally, for not necessarily linear systems, if there exists some
nonnegative definite smooth function $V$ on states with the property that, for
some $\gamma >0$,
\[
\nabla V(x).f(x,u)  \;\leq  \; -y^2 \,+\, \gamma \, uy 
\]
for all $x\in \R^n$, $u\in \R$, and $y=h(x)$, then the system is OSP.
Indeed, integrating along solutions corresponding to $x(0)=0$, and using that
$V$ is nonnegative definite (so that $V(0)=0$ and $V(x(T))\geq 0$),
one has that
\beqn
0 &\leq & V(x(T))-V(0) \\
&\leq & -\int_0^T y(s)^2\,ds \,+\, \gamma  \int_0^T u(s)y(s)\,ds
\eeqn
and thus $\normT{y}^2 \leq  \gamma  \ipT{u}{y}$ as claimed.
This property can be checked by means of nonlinear versions of the KYP Lemma,
see e.g.~\cite{khalil,vdS}.

Yet another way of stating the estimate~(\ref{OSP}) is in terms of integral
quadratic constraints (IQC's), cf.~\cite{megretski-rantzer}:
one may equivalently write ``$w^TMw\geq 0$'' in $L^2$ for i/o pairs $w=(u,y)'$
and where: 
\[
M \;=\;
\pmatrix{0 & \gamma /2\cr
         \gamma /2 & -1}
\]
The powerful tools for analysis of IQC's, based on LMI's, as developed by
Megretski and Rantzer and others, should thus be useful for the study of secant
gains. (We wish to thank R. Sepulchre for suggesting this reformulation.)

\begin{figure}[h,t]
\begin{center}
\includegraphics[width=6cm]{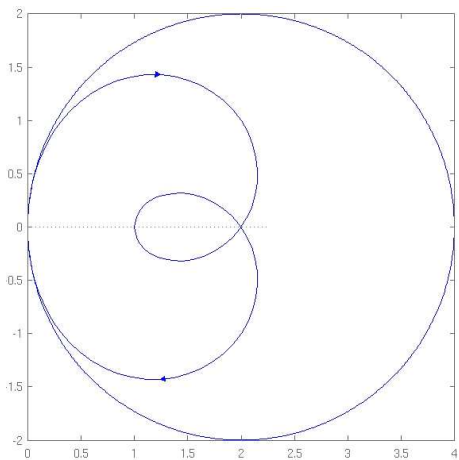}
\caption{Nyquist plot for $\frac{2s+1}{s^2+s+1}$ and circle
$\abs{s-2}\leq 2$, $\gamma _s=4$, $\gi\approx2.25$}
\label{fig-nyquist}
\end{center}
\end{figure}

We pointed out that the induced $L^2$ gain $\gi$ is upper bounded by the
secant gain $\gamma _s$.  In general, one has the strict inequality
$\gi<\gamma _s$.
For example consider the linear system with transfer function
\[
G(s) \;=\; \frac{2s+1}{s^2+s+1} \,.
\]
This is a scalar multiple of $(s+1/2)/(s^2+s+1)$, so it is SPR by the criterion
mentioned earlier, and hence OSP.  Explicitly:
\[
\gamma _s\;=\;
\sup_{\omega \in \R}\frac{\abs{p(i\omega )}^2}{\RE [p(i\omega )\overline{q(i\omega )}]} \;=\;
\sup_{\omega \in \R}\frac{1+4\omega ^2}{1+\omega ^2} \;=\; 4
\]
and
\beqn
\gi&=&\sup_{\omega \in \R}\abs{\frac{1+2i\omega }{1-i\omega -\omega ^2}} \;=\;
        \sup_{\omega \in \R}\sqrt{\frac{1+4\omega ^2}{1-\omega ^2+\omega ^4}}\\
&=& \sqrt{2+(2/3)\sqrt{21}}\;\approx\; 2.25 \;< \;4
\eeqn
(the maximum value is achieved at $\omega =1/2\,\sqrt {\sqrt {21}-1}$).
Graphically, we can see these conclusions from Figure~\ref{fig-nyquist},
 which
shows that the smallest circle of the form $\abs{s-\gamma /2}\leq \gamma /2$ which contains
the Nyquist plot must have $\gamma =4$ (circle shown), so that this is the value of
$\gamma _s$, but $\gi\approx2.25$ because the plot fits in a circle (not shown)
centered at the origin with radius $\approx2.25$.

To conclude, let us provide a direct proof of the main theorem in the linear
case.  This proof, when specialized to linear one-dimensional systems, is
basically the same as the proof given in~\cite{Thron}.
We assume a unity negative feedback about the cascade in
Figure~\ref{fig-system}, where each system has transfer function $G_i(s)$ and
secant gain $\gamma _i$.  As remarked earlier, this means that an estimate
$\abs{G_i(\lambda )}\leq \gamma _i\cos \theta (G_i(\lambda ))$ as in~(\ref{OSP-linear-angle}) holds for
every $\lambda $ with real part $\geq 0$.  If the closed-loop were not to be stable,
then there would exist a pole $\lambda $ with real part $\geq 0$ for $G/(1+G)$, where
$G=G_1\ldots G_n$ is the open-loop system.  For any such $\lambda $:
\be{linear-cl}
G_1(\lambda )G_2(\lambda )\ldots G_n(\lambda ) = -1
\ee
from which we conclude, writing $\theta _i=\theta (G_i(\lambda ))$,
that $\sum_i \theta _i$ is a multiple of $\pi $.
Moreover, taking absolute values in~(\ref{linear-cl}) and using 
$\abs{G_i(\lambda )}\leq \gamma _i\cos \theta _i$, we have also that
\[
1 \;\leq \; \gamma _1\gamma _2\ldots \gamma _n \, \Pi _{i=1}^n \cos \theta _i
\]
and so $\gamma _1\ldots \gamma _n\geq \left(\sec {\pi \over n}\right)^n$, again using the
convexity of $-\ln\cos x$.
Thus no such poles can exist, if the hypothesis of the theorem holds.

\end{document}